# PADDED POLYNOMIALS, THEIR COUSINS, AND GEOMETRIC COMPLEXITY THEORY

HARLAN KADISH AND J.M. LANDSBERG

ABSTRACT. We establish basic facts about the varieties of homogeneous polynomials divisible by powers of linear forms, and explain consequences for geometric complexity theory. This includes quadratic set-theoretic equations, a description of the ideal in terms of the kernel of a linear map that generalizes the Foulkes-Howe map, and an explicit description of the coordinate ring of the normalization. We also prove asymptotic injectivity of the Foulkes-Howe map.

1. INTRODUCTION

**Motivation.** In the *geometric complexity theory* (GCT) proposed by Mulmuley and Sohoni, see [9, 10, 2], one is concerned with determining when one orbit closure contains another. The flagship such problem is as follows: Let $\det_n \in S^n\mathbb{C}^{n^2}$ denote the determinant (a homogeneous polynomial of degree $n$ in $n^2$ variables), and let $\perm_m \in S^m\mathbb{C}^{m^2}$ denote the permanent (a homogeneous polynomial of degree $m$ in $m^2$ variables). Assume $n > m$. Classically one was interested in the value $n = n(m)$ for which $\perm_m$ could be realized a specialization of $\det_n$. In order to use the tools of algebraic geometry and representation theory, one would like the polynomials to live in the same ambient space. To do this, let $\mathbb{C}$ have a linear coordinate $\ell$, and consider any linear inclusion $\mathbb{C}^{m^2} \oplus \mathbb{C} \subset \mathbb{C}^{n^2}$. Then the original question is whether or not $\ell^{n-m} \perm_m$ is in $\End(\mathbb{C}^{n^2}) \cdot \det_n$. Since the endomorphisms do not form a group, Mulmuley and Sohoni proposed to instead consider the question of whether or not $\ell^{n-m} \perm_m \in \overline{GL_{n^2} \cdot \det_n}$, or equivalently, whether or not $\overline{GL_{n^2} \cdot \ell^{n-m} \perm_m} \subset \overline{GL_{n^2} \cdot \det_n}$. There are variants of this question with other polynomials.

One approach to separating the orbit closures is to look for modules of polynomials in the ideal of $\overline{GL_{n^2} \cdot \det_n}$ that are not in the ideal $\overline{GL_{n^2} \cdot \ell^{n-m} \perm_m}$. One can look for such polynomials by searching for ways in which the hypersurface $\{\det_n = 0\}$ is pathological, but $\{\ell^{n-m} \perm_m = 0\}$ is not.

We call a polynomial of the form $P = \ell^{n-m}h$, $h \in S^m\mathbb{C}^M$, a *padded polynomial*. (In the GCT literature, they are often called "blasted polynomials".) A key difficulty is that *any* padded polynomial is pathological by most standard measures (e.g., the codimension of singular locus of the hypersurface $\{P = 0\}$ is zero). So an important step in the GCT program is to determine which modules of polynomials on $S^n\mathbb{C}^N$ vanish on padded polynomials. This motivates the following definition from [2]:

Throughout this paper, $W$ denotes a complex vector space of dimension **w**. Define

$F_{n-m}(Sub_k(S^nW^*))$
$\quad := \{P \in S^nW^* \mid \exists W' \subset W^*, \dim W' = k, P = \ell^{n-m}h, \text{for some } \ell \in W', h \in S^mW'\}$,

Kadish supported by DOE ASCR grant DE-SC0002505 (Topology for Statistical Modeling of Petascale Data). Landsberg supported by NSF grant DMS-1006353.

and when $k = \mathbf{w}$, write:
$$F_{n-m}(S^nW^*) := \{P \in S^nW^* \mid P = \ell^{n-m}h, \text{for some } \ell \in W, h \in S^mW\}.$$

Before trying to separate $\det_n$ from $\text{perm}_m$, one might first simply try to separate $\det_n$ from a generic sequence $h_m \in S^m\mathbb{C}^{m^2}$, which motivates the following definition:

**Definition 1.1.** Let $(d_n)$ be a sequence of polynomials $d_n \in S^n\mathbb{C}^{N(n)} =: S^nW_n$ that is complete for some complexity class. Let $f : \mathbb{Z}_+ \to \mathbb{Z}_+$, $M : \mathbb{Z}_+ \to \mathbb{Z}_+$ be increasing functions of $m$. Let $U_n \subset Sym(S^nW_n)$ be a sequence of $GL(W_n)$-modules of polynomials such that $U_n \subset I(\overline{GL_N \cdot d_n})$. We say the sequence of modules $U_n$ is $(f(m), M(m))$-GCT useful if $U_{f(m)} \not\subset I(F_{f(m)-m}(Sub_{M(m)+1}(S^{f(m)}W^*_{f(m)})))$.

In the flagship case, $d_n = \det_n$ and $M(m) = m^2$. For this case, the best known sequence of modules is obtained via the equations for varieties with degenerate dual varieties, which yield $(\frac{m^2}{2}, m^2)$-GCT useful equations; see [6]. In fact, the permanent is generic with respect to these equations, so one obtains a lower bound on the (GCT) determinantal complexity of the permanent. (Note that there is no standard definition for the dual variety of a non-reduced variety.)

Separating $\det_n$ from a generic sequence $h_m \in S^m\mathbb{C}^{m^2}$ would be done by finding a sequence of modules that are $(m^c, m^2)$-GCT useful for all $c > 0$. By a GCT analogue of the Razborov-Rudich *natural proof barrier* [11], this alone is unlikely to be sufficient to separate $\det_n$ from $\text{perm}_m$, but there is little hope of separating the determinant from the permanent before this problem is solved.

**Results.** Recall that polynomial representations of $GL(W)$ are indexed by partitions $\pi$, where the number of parts of $\pi$, denoted $\ell(\pi)$, is at most $\mathbf{w}$. We let $S_\pi W$ denote the corresponding irreducible $GL(W)$-module. We write $S_{\overline{\pi}}W$ to denote a specific realization of $S_\pi W$ in $S^d(S^nW)$.

**Theorem 1.2.** *A necessary condition for a module $S_{\overline{\pi_n}}W_n \subset I_d(\overline{GL_N \cdot d_n})$ to be $(n(m), M(m))$-GCT useful is*

(1) $\ell(\pi_{n(m)}) \leq M(m) + 1$,
(2) *If $\pi_n = (p_1, \ldots, p_t)$, then $p_1 \geq d(n - m)$.*

*Moreover, if $p_1 \geq \min\{d(n-1), dn - m\}$, then the necessary conditions are also sufficient. In particular, for $p_1$ sufficiently large, GCT usefulness depends only on the partition $\pi$, not how the module $S_\pi W$ is realized as a space of polynomials.*

The equations of $F_{n-m}(Sub_k(S^nW^*))$ are easily deduced from those of $F_{n-m}(S^nW^*)$, see Proposition 1.12, which explains condition (1). As a consequence, from now on we only consider the case $k = \mathbf{w}$. Condition (2) and the "moreover" assertion are a consequence of the following theorem:

**Theorem 1.3.** $I_d(F_{n-m}(S^nW^*))$ *contains the isotypic component of $S_\pi W$ in $S^d(S^nW)$ for all $\pi = (p_1, \ldots, p_d)$ with $p_1 < d(n-m)$. It does not contain a copy of any $S_\pi W$ where $p_1 \geq \min\{d(n-1), dn - m\}$.*

Theorem 1.3 is proved in §2.

*Remark* 1.4. For the so-called subspace variety $F_0(Sub_k(S^nW^*))$, a module appears in its ideal if and only if its entire isotypic component in $Sym(S^nW)$ appears. This property fails for $F_{n-m}(S^nW^*)$. For example, consider $I_3(F_{n-2}(S^n\mathbb{C}^2))$. For $n \geq 6$, the module $S_{(3(n-2),6)}W$ appears with multiplicity two in $S^3(S^nW)$, by the stability result in [7, Thm 4.2.2]. Using [13, Cor. 4a], one computes that the module $S_{(3(n-2),6)}W$ has multiplicity one in $\mathbb{C}[F_{n-2}(S^n\mathbb{C}^2)]$ and hence multiplicity one in $I_3(F_{n-2}(S^n\mathbb{C}^2))$ as well. Thus as far as a GCT guide for "where to look" for good (abstract) modules is concerned, Theorem 1.2 appears to be sharp.

**Question 1.5.** For any $S_\pi W$ occurring in $S^d(S^n W)$ with $p_1 \geq d(n-m)$, is there always some copy of $S_\pi W$ occurring in the coordinate ring $\mathbb{C}[F_{n-m}(S^n W^*)]$?

Regarding set-theoretic equations, we have the following result:

**Proposition 1.6.** *If $4m < n$, then the variety of padded polynomials $F_{n-m}(S^n W^*)$ is cut out set-theoretically by equations of degree two. For any $n > m$,*

$$I_2(F_{n-m}(S^n W^*)) = \bigoplus_{k=1}^{\lfloor \frac{n}{2} \rfloor + m} S_{2(n-m-k), 2(m+k)} W. \tag{1}$$

As explained below, Proposition 1.6 can be deduced from from [13, Thm. 2.15] and Proposition 1.12. We give a direct proof of the second assertion in §2.

It is useful to study the following more general class of varieties, which interpolate between the Chow variety of polynomials that are products of linear forms, and the variety $F_{n-m}(S^n W^*)$: Let $s_1 \geq \cdots \geq s_q \geq 0$. We use the notations $\mathbf{s} = (s_1, \ldots, s_q)$, $(s_1, \ldots, s_q) = ((q_1)^{t_1}, \ldots, (q_f)^{t_f})$ where $s_1 = \cdots s_{t_1} = q_1$, etc., and $|\mathbf{s}| = s_1 + \cdots + s_q$ Consider the varieties

$$F_\mathbf{s}(S^n W^*) = F_{s_1,\ldots,s_q}(S^n W^*) := \{P \in S^n W^* \mid P = \ell_1^{s_1} \cdots \ell_q^{s_q} h, \ \ell_j \in W^*, \ h \in S^{n-(s_1+\cdots+s_q)} W^*\}.$$

The Chow variety of products of linear forms is $F_{(1)^n}(S^n W^*)$.

**Theorem 1.7.** *For all $\delta > 1$, there is a $GL(W)$-module map*

$$\mathcal{F}_\mathbf{s}^\delta : S^\delta(S^n W) \to S^{t_1}(S^{\delta q_1} W) \otimes \cdots \otimes S^{t_f}(S^{\delta q_f} W) \otimes S^\delta(S^{n-|\mathbf{s}|} W) \tag{2}$$

*whose kernel is $I_\delta(F_\mathbf{s}(S^n W^*))$. The map $\mathcal{F}_\mathbf{s}^\delta$ is given explicitly in §3.*

The map $\mathcal{F}_\mathbf{s}^d$ is a generalization of the Foulkes-Howe map $S^\delta(S^n W) \to S^n(S^\delta W)$ (case $\mathbf{s} = (1)^n$) whose kernel is the ideal in degree $\delta$ of the Chow variety (see, e.g., [5, §8.6] or [4, Chap. 4]), and the map $\theta$ of [13], which is the case $\mathbf{w} = 2$ and $\mathbf{s} = s_1$. Theorem 1.7 is proved in §3.

The following proposition generalizes results of [2, §6.4] which treats the case $q = 1$. It is proved in §4.

**Proposition 1.8.** *Assume $s_1 > s_2 > \cdots > s_q$. Suppose $\dim W \geq 3$. Let $Nor(\mathbb{C}[F_\mathbf{s}(S^n W^*)])$ denote the normalization of $\mathbb{C}[F_\mathbf{s}(S^n W^*)]$ and note that it is a graded ring. Then*

$$Nor(\mathbb{C}[F_\mathbf{s}(S^n W^*)])_\delta = S^{\delta s_1} W \otimes \cdots \otimes S^{\delta s_q} W \otimes S^\delta(S^{n-\mathbf{s}} W).$$

*Remark* 1.9. Note that in the case $s_1 > s_2 > \cdots > s_q$, the target of the map (2) is $Nor(\mathbb{C}[F_\mathbf{s}(S^n W^*)])_\delta$.

Thus the situation for padded polynomials, at least in the range relevant for GCT, is similar to that of the Chow variety - one has a description of the ideal in any degree in terms of the kernel of an explicit linear map, and one has an explicit construction of a collection set-theoretic equations. In fact the situation is slightly better: one has the set-theoretic equations described explicitly as $GL(W)$-modules, this collection of modules is multiplicity free, and there is a very simple description of a large part of the ideal. (In the Chow case it is known that the modules are multiplicity-free, but the modules are not known explicitly.)

For a linear map $f : A \to B$, let $f^T : B^* \to A^*$ denote its transpose. We write $\mathcal{F}_d^\delta = {}^A\mathcal{F}_d^\delta : S^\delta(S^d A) \to S^d(S^\delta A)$ to emphasize the vector space $A$.

**Proposition 1.10.** *The Foulkes-Howe map* $\mathcal{F}^\delta_{(1)^d} : S^\delta(S^dW) \to S^d(S^\delta W)$ *is self-dual, in the sense that for* $Q \in S^d(S^\delta W^*)$,

$$({}^W\mathcal{F}^\delta_{(1)^d})^T(Q) = {}^{W^*}\mathcal{F}^d_{(1)^\delta}(Q).$$

Proposition 1.10 is proved in §5.

A conjecture that has come to be called the "Foulkes-Howe conjecture" states that the map $\mathcal{F}^\delta_{(1)^\delta}$ is injective for $\delta \leq d$ and surjective for $\delta \geq d$. The conjecture is known to be false as stated when $\mathbf{w} = d = \delta = 5$ as a consequence of results in [8], although the kernel in that case is not known explicitly. It is also known that the map is eventually surjective. More precisely, [1, Thm. 3.3] states that when

$$(3) \qquad d \geq (\delta - 1)(\mathbf{w} - 1)((\delta - 1)\lceil \frac{\binom{\delta + \mathbf{w} - 1}{\mathbf{w} - 1}}{\mathbf{w}} \rceil + \delta)$$

the map $\mathcal{F}^d_{(1)^\delta}$ is surjective. We immediately obtain, as a Corollary of Proposition 1.10:

**Corollary 1.11.** *When (3) holds, the map* $\mathcal{F}^\delta_{(1)^d} : S^\delta(S^dW) \to S^d(S^\delta W)$ *is injective.*

**History and first steps towards the proofs.** The variety $F_{n-m}(S^n\mathbb{C}^2)$ was an object of classical study, as it corresponds to polynomials with multiple roots. In a series of papers [12, 13, 14], J. Weyman determined extensive information about its ideal, Hilbert function, etc. In particular, he showed that when $4m < n$, the ideal of $F_{n-m}(S^n\mathbb{C}^2)$ is generated in degree two [13, Thm. 2.15], and he gave a recipe for determining the ideal in any degree [13, Cor. 4a]. Now, a polynomial $P$ on $W$ divides a polynomial $Q$ iff $P$ restricted to every $\mathbb{C}^2 \subset W$ divides $Q$ restricted to the $\mathbb{C}^2$. This remark combined with [13, Thm. 2.15] proves the first assertion of Proposition 1.6. The varieties $F_\mathbf{s}(S^nW^*)$ are discussed in [3] where a recipe for set-theoretic equations for all cases is given. The equations there are in terms of covariants. In particular, they are of degree greater than two for $F_{n-m}(S^nW^*)$.

**Inheritance.** For a partition $\pi = (p_1, \ldots, p_\mathbf{w})$, we write $|\pi| = p_1 + \cdots + p_\mathbf{w}$ and observe that $S_\pi W$ appears in $W^{\otimes|\pi|}$ with multiplicity $\dim[\pi]$, where $[\pi]$ is the corresponding representation of $\mathfrak{S}_{|\pi|}$, the symmetric group on $|\pi|$ elements. Let $M_\pi \subset W^{\otimes|\pi|}$ denote the isotypic component associated to $S_\pi W$, and use the notation $S_{\overline{\pi}} W \subset M_\pi$ for a particular copy of $S_\pi W$, given by some linear combination of Young symmetrizers. Note that if $k < \mathbf{w}$ and we have an inclusion $\mathbb{C}^k \subset W$, then the module $S_{\overline{\pi}}\mathbb{C}^k$ induces a unique module $S_{\overline{\pi}}W$. We recall the following proposition:

**Proposition 1.12.** *(Inheritance, see [5, §7.4])* Notations as above. $I_d(F_{n-m}(Sub_k(S^nW^*)))$ consists of all modules $S_{\overline{\pi}}W$ such that $S_{\overline{\pi}}\mathbb{C}^k$ is in the ideal of $F_{n-m}(S^n\mathbb{C}^{k^*}))$ and all modules whose associated partition has length at least $k + 1$.

Now, $S^k(S^nW)$ only consists of modules $S_\pi W$ with $\ell(\pi) \leq k$. So when studying the ideal in degree $k$ of a subvariety of $S^nW$, by the remarks above we may restrict our study to its cousin in $S^n\mathbb{C}^k$.

**Acknowledgments.** We thank Michel Brion, Jaydeep Chipalkatti, and Jerzy Weyman for useful conversations.

## 2. Proof of Theorems 1.3, 1.2 and equation (1)

Since all partitions in $S^d(S^nW)$ have length at most $d$, we may assume by Proposition 1.12 that $\mathbf{w} = d$. Fix a (weight) basis $e^1, \ldots, e^d$ of $W^*$ with dual basis $x_1, \ldots, x_d$ of $W$. Note any element $\ell^{n-m}h \in F_{n-m}(S^nW^*)$ is in the $GL(W)$-orbit of $(e^1)^{n-m}h$, so it will be sufficient to show that the ideal in degree $d$ contains the modules vanishing on the orbits of elements of the form $(e^1)^{n-m}h$. The highest weight vector of any copy of $S_{(p_1,\ldots,p_d)}W$ in $S^d(S^nW)$ will be a linear combination of vectors of the form $m_I := (x_1^{i_1^1} \cdots x_d^{i_d^1}) \circ \cdots \circ (x_1^{i_1^d} \cdots x_d^{i_d^d})$, where $i_j^1 + \cdots + i_j^d = p_j$ and $i_1^j + \cdots + i_d^j = n$ for all $j$. All the $m_I$ vanish on any $(e^1)^{n-m}h$ unless $p_1 \geq d(n-m)$, proving $I_d(F_{n-m}(S^nW))$ contains all modules $S_\pi W \subset S^d(S^nW)$ with $p_1 < d(n-m)$. This proves the first assertion of Theorem 1.3.

For a coordinate-free proof, apply the Pieri rule to the the target of the $GL(W)$-module map

$$\mathcal{F}_{n-m}^d : S^d(S^nW) \to S^{d(n-m)}W \otimes S^d(S^mW)$$

of Theorem 1.7, which is independent of Theorem 1.3. If $S_\nu W$ is a submodule of the target, then there is a submodule $S_\mu W \subset S^d(S^mW)$ such that the Young diagram of $\nu = (\nu_1, \ldots, \nu_d)$ can be formed by adding $d(n-m)$ boxes to the diagram of $\mu$, with no two of these new boxes in the same column. Hence the Young diagram of $\nu$ must have at least $d(n-m)$ columns, that is, $\nu_1 \geq d(n-m)$. By Schur's Lemma, if $S_\pi W \subset S^d(S^nW)$ and $\pi = (p_1, \ldots, p_d)$ has $p_1 < d(n-m)$, then $S_\pi W$ lies in the kernel of $\mathcal{F}_{\mathbf{s}}^d$. □

To see the second assertion of Theorem 1.3, define a map on basis elements by

$$(4) \qquad m_{x_1} : S^d(S^mW) \to S^d(S^nW)$$

$$(x_1^{i_1^1} x_2^{i_2^1} \cdots x_d^{i_d^1}) \circ \cdots \circ (x_1^{i_1^d} \cdots x_d^{i_d^d}) \mapsto (x_1^{i_1^1+(n-m)} x_2^{i_2^1} \cdots x_d^{i_d^1}) \circ \cdots \circ (x_1^{i_1^d+(n-m)} \cdots x_d^{i_d^d})$$

and extend linearly. Note that $m_{x_1}$ takes a vector of weight $\mu = (q_1, q_2, \ldots, q_d)$ to a vector of weight $\pi = (p_1, \ldots, p_d) := \mu + (d(n-m)) = (q_1 + d(n-m), q_2, \ldots, q_d)$ in $S^d(S^nW)$. Moreover, if that vector is a highest weight vector, its image is also a highest weight vector, as raising operators annihilate $x_1$. The image of the space of highest weight vectors for $\mu$ in $S^d(S^mW)$ under $m_{x_1}$ lies in the space of highest weight vectors for $\mu + (d(n-m))$ in $\mathbb{C}[F_{n-m}(S^nW^*)]$ because, for example, such a polynomial will not vanish on $(e^1)^{n-m}[(e^1)^{i_1^1} \cdots (e^d)^{i_d^1} + \cdots + (e^1)^{i_1^d} \cdots (e^d)^{i_d^d}]$.

By [7, Thm. p1],

$$\mathrm{mult}(S_\pi W, S^d(S^mW)) = \mathrm{mult}(S_{\pi+d(n-m)}W, S^d(S^nW))$$

as soon as $d \geq q_2 + \cdots + q_d$. This condition implies $q_1 \geq dm - d$, so $p_1 = q_1 + dn - dm \geq dn - d = d(n-1)$ as required. For the sufficiency of $p_1 \geq dn - m$, note that if $p_1 \geq (d-1)n + (n-m) = dn - m$, then in a tensor of weight $\pi$, each of the exponents $i_1^1, \ldots, i_1^d$ of $x_1$ must be at least $n - m$. So there again exists a polynomial in $F_{n-m}(S^nW)$ where a vector of weight $\pi$ does not vanish. □

To prove equation (1), write $x = x_1, y = x_2$, and note that the highest weight vector of $S_{2n-2s,2s}W \subset S^2(S^nW)$ will be a linear combination of vectors of the form

$$(x^n) \circ (x^{n-2s}y^{2s}), (x^{n-1}y^1) \circ (x^{n-2s+1}y^{2s-1}), \ldots, (x^{n-s}y^s) \circ (x^{n-s}y^s)$$

with a nonzero coefficient on the last factor. These vectors vanish on $e_1^{n-m}h$ for any $h$ if and only if $s \geq m + 1$. □

## 3. Proof of Theorem 1.7

Define $\mathcal{F} = \mathcal{F}_{\mathbf{s}}^\delta$ as follows: Include $S^\delta(S^dW) \subset W^{\delta d}$ and label the entries of a tensor as

$$(W_{1,1} \otimes \cdots \otimes W_{1,d}) \otimes (W_{2,1} \otimes \cdots \otimes W_{2,d}) \otimes \cdots \otimes (W_{\delta,1} \otimes \cdots \otimes W_{\delta,d}),$$

where all $W_{i,j}$ are copies of $W$. Extract the first $s_1$ entries of each group in parentheses and collect them in the beginning of a new expression, followed by entries $s_1 + 1, \ldots, s_1 + s_2$, and so on through $s_q$:

$$(W_{1,1} \otimes \cdots \otimes W_{1,s_1} \otimes \cdots \otimes W_{\delta,1} \otimes \cdots \otimes W_{\delta,s_1})$$
$$\otimes (W_{1,s_1+1} \otimes \cdots \otimes W_{1,s_1+s_2} \otimes \cdots \otimes W_{\delta,s_1+1} \otimes \cdots \otimes W_{\delta,s_1+s_2})$$
$$\otimes \cdots \otimes (W_{1,|\mathbf{s}|+1} \otimes \cdots \otimes W_{1,d} \otimes \cdots \otimes W_{\delta,|\mathbf{s}|+1} \otimes \cdots \otimes W_{\delta,d})$$

Next, symmetrize each group of $\delta s_i$ and each group of $d - |\mathbf{s}|$ entries to obtain an element of $S^{\delta s_1} \otimes \cdots \otimes S^{\delta s_q}W \otimes (S^{d-\mathbf{s}}W)^{\otimes \delta}$. Symmetrize the ordering of the $\delta$ groups of $d - |\mathbf{s}|$ entries to obtain an element in $S^{\delta s_1}W \otimes \cdots \otimes S^{\delta s_q}W \otimes S^\delta(S^{d-|\mathbf{s}|}W)$. Finally, for each group of equal $s_j$'s, symmetrize the groups to get an element of $S^{t_1}(S^{\delta q_1}W) \otimes \cdots \otimes S^{t_f}(S^{\delta q_f}W) \otimes S^\delta(S^{d-|\mathbf{s}|}W)$.

To prove $\ker \mathcal{F} = I_\delta(F_\mathbf{s}(S^dW^*))$, observe that $S^\delta(S^nW)$ is spanned by terms of the form $P = x_1^n \cdots x_\delta^n$. Hence without loss of generality, we show $P = x_1^d \cdots x_\delta^d \in I_\delta(F_\mathbf{s}(S^dW^*))$ if and only if $P \in \ker \mathcal{F}$. To that end, choose $\alpha_1^{s_1} \cdots \alpha_q^{s_q} h \in F_{s_1,\ldots,s_q}(S^dW)$, where $\alpha_i \in W$ and $h \in S^{d-|\mathbf{s}|}W$. Letting $\overline{P}$ be the polarization of $P$,

$$\begin{aligned}
\overline{P}(\alpha_1^{s_1} \cdots \alpha_q^{s_q} h) &= x_1^d(\alpha_1^{s_1} \cdots \alpha_q^{s_q} h) \cdots x_\delta^d(\alpha_1^{s_1} \cdots \alpha_q^{s_q} h) \\
&= x_1^{|\mathbf{s}|}(\alpha_1^{s_1} \cdots \alpha_q^{s_q}) x_1^{d-|\mathbf{s}|}(h) \cdots x_\delta^{|\mathbf{s}|}(\alpha_1^{s_1} \cdots \alpha_q^{s_q}) x_\delta^{d-|\mathbf{s}|}(h) \\
&= \left(\prod_{i,j} x_i(\alpha_j)^{s_j}\right) \cdot x_1^{d-|\mathbf{s}|}(h) \cdots x_\delta^{d-|\mathbf{s}|}(h)
\end{aligned}$$

On the other hand,

$$\mathcal{F}(P) = \overbrace{x_1 \cdots x_1}^{s_1 \text{ terms}} \cdots \overbrace{x_\delta \cdots x_\delta}^{s_1 \text{ terms}} \otimes \cdots \otimes \overbrace{x_1 \cdots x_1}^{s_q \text{ terms}} \cdots \overbrace{x_\delta \cdots x_\delta}^{s_q \text{ terms}} \otimes (x_1^{d-|\mathbf{s}|}) \cdots (x_\delta^{d-|\mathbf{s}|})$$
$$\in S^{\delta s_1} \otimes \cdots \otimes S^{\delta s_q}W \otimes S^\delta(S^{d-|\mathbf{s}|}W),$$

hence $\overline{P}(\alpha_1^{s_1} \cdots \alpha_q^{s_q} h) = \langle \mathcal{F}(P), \alpha_1^{\delta s_1} \otimes \cdots \otimes \alpha_q^{\delta s_q} \otimes h^\delta \rangle$. Thus if $P \in \ker \mathcal{F}$, then $P$ vanishes on $F_\mathbf{s}(S^dW)$. If $\mathcal{F}(P) \neq 0$, then there exist $\alpha_i \in W, h \in S^{d-|\mathbf{s}|}W$ such that the pairing $\langle \mathcal{F}(P), \alpha_1^{\delta s_1} \otimes \cdots \otimes \alpha_q^{\delta s_q} \otimes h^\delta \rangle$ is nonzero. The result follows. □

## 4. Proof of Proposition 1.8

Note first that if $d - |\mathbf{s}| = 1$, then $F_{s_1,\ldots,s_q}(S^dW) = F_{s_1,\ldots,s_q,1}(S^dW)$. So without loss of generality we assume $d - |\mathbf{s}| \geq 2$.

We construct a Kempf-Weyman desingularization (see [15, Chap. 5]). Let $\mathcal{S} \to \mathbb{P}W$ denote the tautological line bundle corresponding to $\mathcal{O}(-1)_{\mathbb{P}W}$. Let $E := \mathcal{S}^{s_1} \boxtimes \cdots \boxtimes \mathcal{S}^{s_q} \otimes S^{d-|\mathbf{s}|}W \to (\mathbb{P}W)^{\times q}$, where $\boxtimes$ denotes exterior tensor product. This $E$ is a sub-bundle of the trivial bundle over $(\mathbb{P}W)^{\times q}$ with fiber $V = S^{s_1}W \otimes \cdots \otimes S^{s_q}W \otimes S^{d-|\mathbf{s}|}W$, namely,

$$E = \left\{ (\ell_1^{s_1} \otimes \cdots \otimes \ell_q^{s_q} \otimes h, ([\ell_1], \ldots, [\ell_q])) \in V \times (\mathbb{P}W)^{\times q} \right\}.$$

Let $q \colon E \to S^{s_1}W \otimes \cdots \otimes S^{s_q}W \otimes S^{d-|\mathbf{s}|}W$ be the projection map from the fibers of $E$, and write $q(E) = F'$. Then $q$ is a birational map to $F'$, and $E$ is a desingularization of $F'$. The higher cohomology of $Sym(E^*)$ vanishes, and so by [15, Thm 5.1.3(a)],

$$\begin{aligned}
\mathbb{C}[Nor(F')]_\delta &= H^0((\mathbb{P}W)^{\times q}, S^\delta(\mathcal{S}^{s_1*} \otimes \cdots \otimes \mathcal{S}^{s_q*} \otimes S^{n-|\mathbf{s}|}W)) \\
&= S^{\delta s_1}W \otimes \cdots \otimes S^{\delta s_q}W \otimes S^\delta(S^{d-|\mathbf{s}|}W).
\end{aligned}$$

The symmetrization map $\rho\colon V \to S^d W$ provides a surjection map of cones $F' \to F_{\mathbf{s}}(S^d W^*) =: F$. As a projection centered at a subspace disjoint from $F'$, this map is finite. Hence the varieties $F$ and $F'$ will have the same normalization if $\rho$ is generically injective.

Consider the open set $U = F - F_{s_1,\ldots,s_q,1}(S^d W)$ which is dense in $F$ and nonempty because $d - |\mathbf{s}| \geq 2$. If $f \in U$, then $f = \alpha_1^{s_1} \cdots \alpha_q^{s_q} h$ for some $\alpha_i \in W$ and for some $h \in S^{d-|\mathbf{s}|} W$ without linear factors. Let $D \subset E$ be the open set with $\ell_i \neq \ell_j$ for all $i \neq j$; then $q(D)$ contains an open set in $F'$ of tensors $\alpha_1^{s_1} \otimes \cdots \otimes \alpha_q^{s_q} \otimes h$ such that $\alpha_i \neq \alpha_j$ for all $i \neq j$. Since $s_i \neq s_j$ for all $i \neq j$, it follows that $\rho$ is injective on $\rho^{-1}(U) \cap q(D)$ as required. $\square$

## 5. Proof of Proposition 1.10

Let $\langle \cdot, \cdot \rangle : W \times W^* \to \mathbb{C}$ denote the pairing between a vector space and its dual. Choose $Q \in S^d(S^\delta W^*)$. We want to show $({}^W \mathcal{F}_d^\delta)^T(Q) = {}^{W^*}\mathcal{F}_\delta^d(Q)$, that is, for all $P \in S^\delta(S^d W)$, that $\langle {}^W \mathcal{F}_d^\delta(P), Q \rangle = \langle P, ({}^{W^*}\mathcal{F}_\delta^d)(Q) \rangle$. It is sufficient to assume $Q = \alpha_1^\delta \cdots \alpha_d^\delta$ for some $\alpha_j \in W^*$ because such elements span $S^d(S^\delta W^*)$. Similarly, we may assume $P = x_1^d \cdots x_\delta^d$.

$$\begin{aligned}
\langle {}^W \mathcal{F}_d^\delta(P), Q \rangle &= \langle (x_1 \cdots x_\delta)^d, \alpha_1^\delta \cdots \alpha_d^\delta \rangle \\
&= \alpha_1^\delta(x_1 \cdots x_\delta) \cdots \alpha_d^\delta(x_1 \cdots x_\delta) \\
&= \alpha_1(x_1) \cdots \alpha_1(x_\delta) \cdots \alpha_d(x_1) \cdots \alpha_d(x_\delta) \\
&= \langle x_1^d, \alpha_1 \cdots \alpha_d \rangle \cdots \langle x_\delta^d, \alpha_1 \cdots \alpha_d \rangle \\
&= \langle x_1^d \cdots x_\delta^d, (\alpha_1 \cdots \alpha_d)^\delta \rangle \\
&= \langle x_1^d \cdots x_\delta^d, {}^{W^*}\mathcal{F}_\delta^d(Q) \rangle. \quad \square
\end{aligned}$$

*E-mail address*: `hmkadish@math.tamu.edu, jml@math.tamu.edu`